\documentclass{article}
\usepackage{amsmath}
\usepackage{amssymb}
\usepackage{latexsym}
\usepackage{verbatim}
\usepackage{graphics}
\usepackage{amsfonts}
\usepackage{graphicx}
\usepackage{epstopdf}
\newtheorem{thm}{Theorem}
\newtheorem{conj}{Conjecture}
\newtheorem{prop}{Proposition}
\newtheorem{lem}{Lemma}
\newtheorem{defn}{Definition}
\newtheorem{cor}{Corollary}

\begin{document}

\title{Intrinsic Linking and Knotting in Virtual Spatial Graphs}

\author{Thomas Fleming \thanks{Partially supported by a Sigma-Xi grant-in-aid of research, number G20059161614561837}\\
	Department of Mathematics\\
   	University of California, San Diego\\
   	La Jolla, CA 92093-0112\\
	{\it tfleming@math.ucsd.edu}
   \and
	Blake Mellor\\
	Mathematics Department\\
   	Loyola Marymount University\\
   	Los Angeles, CA  90045-2659\\
	{\it bmellor@lmu.edu}}

\date{}

\maketitle

\begin{abstract}

We introduce a notion of intrinsic linking and knotting for virtual spatial graphs.  Our theory gives two filtrations of the set of all graphs, allowing us to measure, in a sense, how intrinsically linked or knotted a graph is; we show that these filtrations are descending and non-terminating.  We also provide several examples of intrinsically virtually linked and knotted graphs.  As a byproduct, we introduce the {\it virtual unknotting number} of a knot, and show that any knot with non-trivial Jones polynomial has virtual unknotting number at least 2.

\end{abstract}

\tableofcontents

\section{Introduction} \label{S:intro}

A {\it spatial graph} is an embedding of a graph $G$ in $R^3$.  Spatial graphs are a natural generalization of knots, which can be viewed as the particular case when $G$ is an $n$-cycle.  Over the past twenty years, there has been considerable work looking at linked and knotted cycles in spatial graphs, beginning with Conway and Gordon's proof that every embedding of $K_6$ has a pair of linked cycles, and every embedding of $K_7$ has a knotted cycle \cite{cg}.  We say that these are examples of {\it intrinsically linked} and {\it intrinsically knotted} graphs, respectively.  Since Conway and Gordon's paper, Robertson, Seymour and Thomas \cite{rst} have classified the intrinsically linked graphs, but the problem of intrinsically knotted graphs is still open, and is an active area of research, along with various variations on these problems (see, for example, \cite{ffnp, fo1, fo2}).

Another recent approach to generalizing knot theory is Kauffman's theory of {\it virtual knots} \cite{ka}.  In previous work, the authors combined these ideas to introduce {\it virtual spatial graphs} \cite{fm}.  The purpose of this paper is to extend some of the problems and theory of classical spatial graphs to this new realm of virtual spatial graphs.  In particular, we will define a notion of {\it intrinsically virtually linked} graphs of various degrees, show that these induce a descending filtration on the set of graphs, and prove some existence results about the filtration.  This gives us a way to talk about ``degree'' of intrinsic linking - for example, in this theory $K_6$ is ``more'' intrinsically linked than the Petersen graph.

In Section \ref{S:virtualgraph} we review the definition of virtual spatial graphs.  In Section \ref{S:linking} we introduce the notion of an {\it n-intrinsically virtually linked (nIVL)} graph  and give examples of graphs which are $nIVL$ but not $(n+1)IVL$ for every $n \geq 1$.  In Section \ref{S:unknotting} we take a short detour into virtual knots to define the {\it virtual unknotting number} and show that every classical knot with non-trivial Jones polynomial has virtual unknotting number at least 2.  Finally, in Section \ref{S:knotting} we define {\it n-virtually intrinsically knotted (nIVK)} graphs, show that all known intrinsically knotted graphs are also $1IVK$, and give examples of graphs which are $(2n-1)IVK$ but not $(2n)IVK$ for every $n \geq 1$.

\section{Virtual Spatial Graphs} \label{S:virtualgraph}

\begin{figure}
    $$\includegraphics{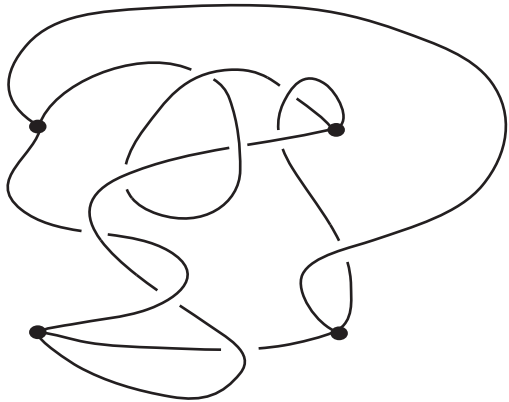}$$
    \caption{A graph diagram} \label{F:diagram}
    \end{figure}

In this section, we will briefly review the definition of a virtual spatial graph, from \cite{fm}.  First, we recall the definition of a classical spatial graph.  A {\it graph} is a pair $G = (V, E)$ of a set of {\it vertices V} and {\it edges} $E \subset V \times V$.  Unless otherwise stated, our graphs are connected and {\it directed}, so that each edge is an {\it ordered} pair of vertices.  An embedding of $G$ in $\mathbb{R}^3$ maps the vertices of $G$ to points in $\mathbb{R}^3$ and an edge $(u, v)$ to an arc in $\mathbb{R}^3$ whose endpoints are the images of the vertices $u$ and $v$, and that is oriented from $u$ to $v$.  We will consider these embeddings modulo equivalence by ambient isotopy.  We can always represent such an embedding by projecting it to a plane so that each vertex neighborhood is a collection of rays with one end at the vertex and crossings of edges of the graph are transverse double points in the interior of the edges (as in the usual knot and link diagrams) \cite{ka2}.  An example of such a diagram is shown in Figure~\ref{F:diagram}.

\begin{figure} 
    $$\includegraphics{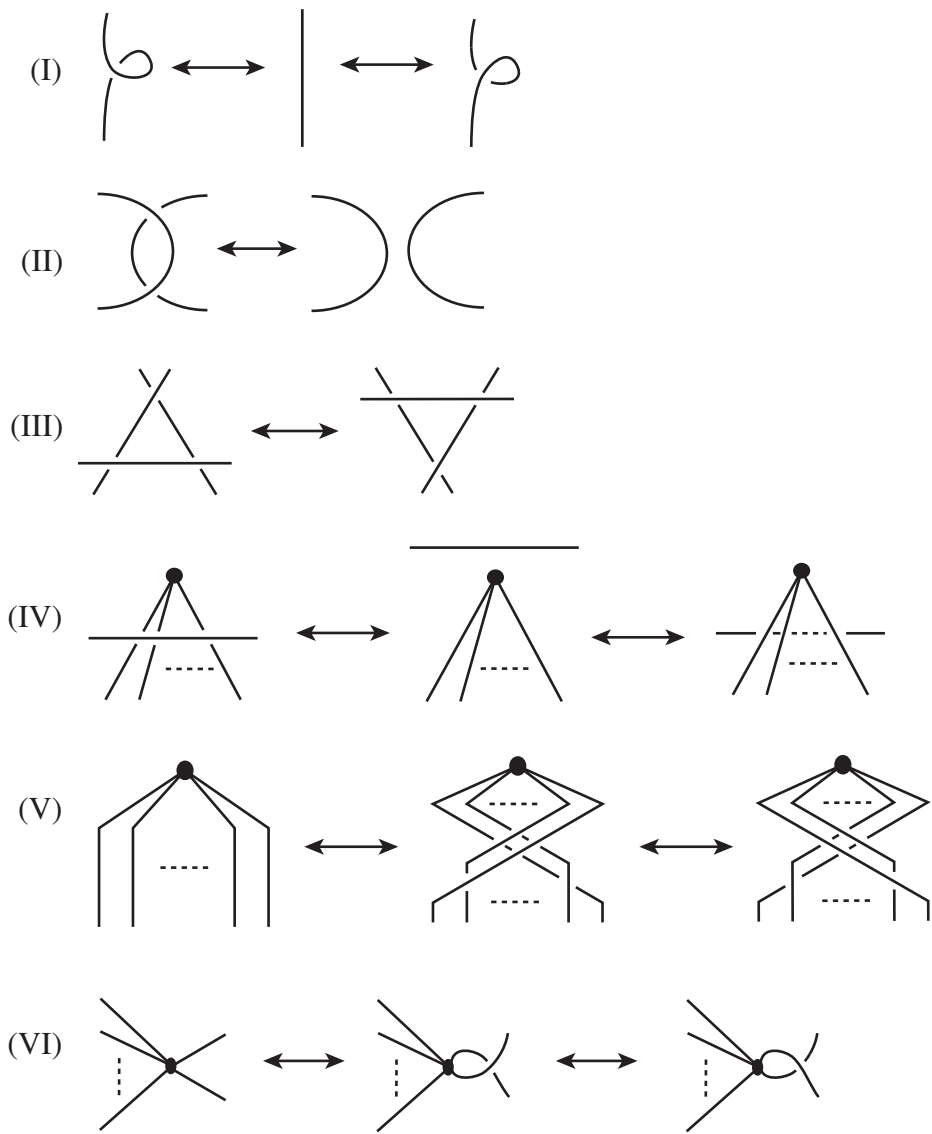}$$
    \caption{Reidemeister moves for graphs} \label{F:classicalmoves}
    \end{figure}

Kauffman \cite{ka2} and Yamada \cite{ya} have shown that ambient isotopy of spatial graphs is generated by a set of local moves on these diagrams which generalize the Reidemeister moves for knots and links.  These Reidemeister moves for graphs are shown in Figure~\ref{F:classicalmoves}.  The first five moves (moves (I) - (V)) generate {\it rigid vertex isotopy}, where the cyclic order of the edges around each vertex is fixed.  Moves (I) - (VI) generate {\it pliable isotopy}, where the order of the vertices around each edge can be changed using move (VI).

\begin{figure}
    $$\includegraphics{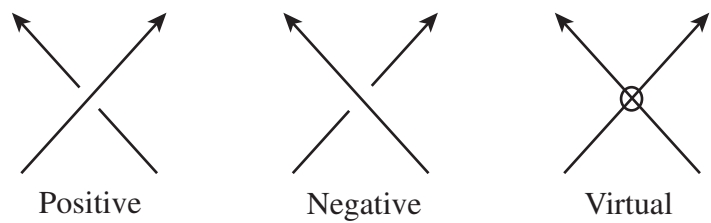}$$
    \caption{Types of crossings} \label{F:crossings}
    \end{figure}

A {\it virtual graph diagram} is just like a classical graph diagram, with the addition of {\it virtual crossings}.  We will represent a virtual crossing as an intersection of two edges surrounded by a circle, with no under/over information.  So we now have three kinds of crossings:  positive and negative classical crossings and virtual crossings (see Figure~\ref{F:crossings}).

\begin{figure}
    $$\includegraphics{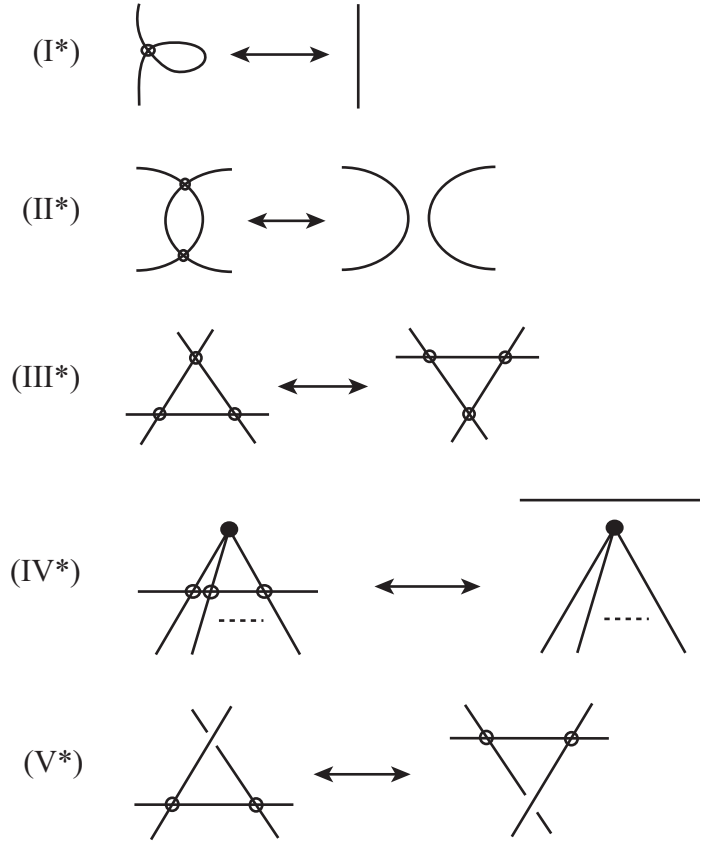}$$
    \caption{Reidemeister moves for virtual graphs} \label{F:virtualmoves}
    \end{figure}

The idea is that the virtual crossings are not really there (hence the name ``virtual").  To make sense of this, we extend our set of Reidemeister moves for graphs to include moves with virtual 
crossings.  We need to introduce 5 more moves, (I*) - (V*), shown in Figure~\ref{F:virtualmoves}.  Notice that moves (I*) - (IV*) are just the purely virtual versions of moves (I) - (IV); move (V*) is the only move 
which combines classical and virtual crossings (in fact, there are two versions of the move, since the classical crossing may be either positive or negative).  We do {\it not} allow the purely virtual version of move (V), since it changes the cyclic order the edges at a vertex (see \cite{fm} for a more detailed discussion).

\section{Intrinsically Virtually Linked Graphs} \label{S:linking}

\begin{figure} 
    $$\includegraphics{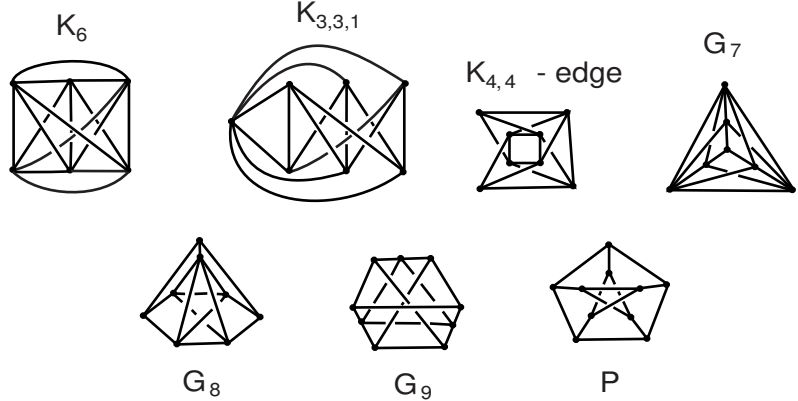}$$
    \caption{The Petersen family of graphs} \label{F:petersenfamily}
    \end{figure}

A graph $G$ is {\it intrinsically linked} if every classical diagram of the graph contains a non-split link whose components are disjoint cycles in the graph.  Robertson, Seymour and Thomas \cite{rst} have shown that every intrinsically linked graph has a minor homeomorphic to a graph in the {\it Petersen family} of graphs (see Figure \ref{F:petersenfamily}).

We will look at {\it virtual} diagrams of the graph $G$.

\begin{defn} \label{D:nIVL}
A graph G is {\bf intrinsically virtually linked of degree n (nIVL)} if every virtual diagram of G with at most n virtual crossings contains a non-trivial virtual link whose components are disjoint cycles in G.
\end{defn}

If we let $IVL_n$ denote the set of intrinsically virtually linked graphs of degree $n$, it is clear that $IVL_0$ is simply the set of classically intrinsically linked graphs, and that $IVL_{n+1}$ is a subset of $IVL_n$.  We would like to know whether it is a {\it proper} subset.  We will find that $IVL_0 = IVL_1$, but that for $n \geq 1$, $IVL_{n+1}$ is a proper subset of $IVL_n$.

We first make the following observation.

\begin{lem} \label{L:odd}
If G is intrinsically linked, then every classical diagram of G contains a pair of linked cycles with odd linking number.
\end{lem}

\noindent{\sc Proof:}  It is known that every diagram of $K_6$ and $K_{3,3,1}$ contains a pair of linked cycles with odd linking number \cite{cg, mo}.  The other graphs in the Petersen family are derived from $K_6$ and $K_{3,3,1}$ by triangle-Y exchanges, in which a triangle (a cycle in the graph of length 3) is replaced by a Y (the three edges of the triangle are removed, and a new vertex is added adjacent to the three vertices of the triangle).  Motwani et. al. \cite{mo} showed that these exchanges preserve the property that every diagram contains a pair of linked cycles with odd linking number, so every graph in the Petersen family has this property.  Since $G$ is intrinsically linked, it has a subgraph $H$ which is homeomorphic to one of the graphs in the Petersen family, by Robertson, Seymour and Thomas \cite{rst}.  So every diagram of $H$ contains a pair of linked cycles with odd linking number, and hence every diagram of G does as well.  $\Box$

\begin{thm} \label{T:0IVL=1IVL}
If G is intrinsically linked, then G is also intrinsically virtually linked of degree 1, so $IVL_0 = IVL_1$
\end{thm}

\noindent {\sc Proof:}  Since $G$ is intrinsically linked, a diagram of $G$ with no virtual crossings will contain a link.  So consider a diagram $V$ of $G$ with exactly one virtual crossing $c$.  Let $V_+$ be the result of replacing $c$ by a positive classical crossing.  By Lemma \ref{L:odd}, $V_+$ contains a pair of linked cycles $C_1$ and $C_2$ with odd linking number.

If the crossing $c$ does not involve an edge of $C_1$ and an edge of $C_2$, then the linking number of $C_1$ and $C_2$ is the same in both $V_+$ and $V$, so $V$ is linked.

On the other hand, if $c$ {\it does} involve an edge of $C_1$ and an edge of $C_2$, then the virtual linking number of $C_1$ with $C_2$ in $V$ will be a half-integer, so $V$ is again linked.

We conclude that every diagram of $G$ with at most one virtual crossing contains a non-trivial link, and hence $G \in IVL_1$. $\Box$

\begin{figure} [h]
    $$\includegraphics{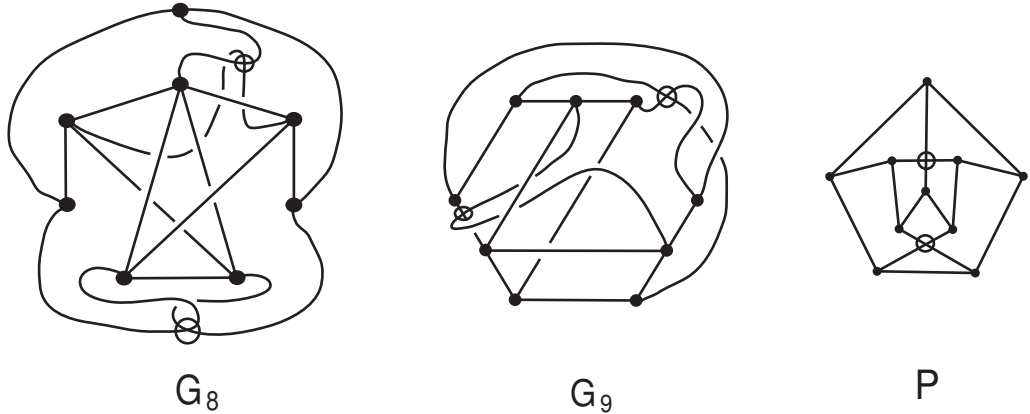}$$
    \caption{$G_8$, $G_9$ and the Petersen graph are not 2IVL} \label{F:magic}
    \end{figure}

Now that we know that $IVL_0 = IVL_1$, is $IVL_2$ any different?  The answer is ``Yes."  Of the graphs in the Petersen family, $G_8$, $G_9$ and the Petersen graph $P$ are in $IVL_1\backslash IVL_2$, while the remaining four are in $IVL_2\backslash IVL_3$.  Figure~\ref{F:magic} shows unlinked virtual diagrams of $G_8$, $G_9$ and $P$ with two virtual crossings, showing that these graphs are not in $IVL_2$.

\begin{figure} [h]
    $$\includegraphics{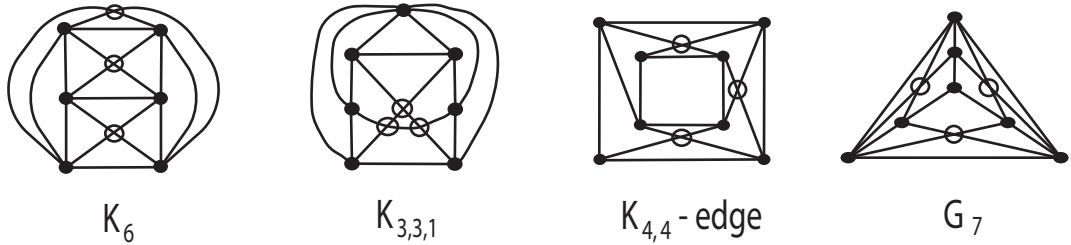}$$
    \caption{$K_6$, $K_{3, 3, 1}$, $K_{4,4}-{\rm edge}$ and $G_7$ are not $3IVL$} \label{F:not3IVL}
    \end{figure}

The diagram for the Petersen graph $P$ illustrates an important observation:  if a graph has minimal crossing number $n$, then it is not in $IVL_n$.  This is because we can simply virtualize every crossing in a minimal diagram, so every link will be purely virtual, and hence trivial.  The remaining graphs of the Petersen family ($K_6$, $K_{3, 3, 1}$, $K_{4,4}-{\rm edge}$ and $G_7$) all have minimal crossing number 3, and so have unlinked diagrams with three virtual crossings, as shown in Figure~\ref{F:not3IVL}.

It only remains to show that these four graphs {\it are} $2IVL$.  We begin with $K_6$.

\begin{prop} \label{P:K6}
$K_6$ is in $IVL_2\backslash IVL_3$.
\end{prop}

\noindent{\sc Proof:}  Assume we have a diagram $G$ of $K_6$ with exactly two virtual crossings $c$ and $d$.  Let $G^*$ be the result of replacing both $c$ and $d$ with classical crossings.  Then $G^*$ contains a pair of disjoint cycles $C_1$ and $C_2$ with odd linking number.  By the argument in Theorem \ref{T:0IVL=1IVL}, $G$ is linked unless both $c$ and $d$ involve an edge of $C_1$ and an edge of $C_2$.  So we can assume that both crossings are between non-adjacent pairs of edges (since $C_1$ and $C_2$ are disjoint cycles).  

We will number the vertices of $K_6$ from 1 to 6.  Without loss of generality, we may assume that crossing $c$ is between edges 1-2 and 4-5.  If crossing $d$ does not occur between the same two edges, then either the link composed of triangles (1-2-3, 4-5-6) or the link composed of triangles (1-2-6, 3-4-5) contains only one virtual crossing between the components (since the two links share only the edges 1-2 and 4-5), and hence has a half-integer virtual linking number in $G$, as in Theorem \ref{T:0IVL=1IVL}.

So we need only consider the case when both $c$ and $d$ are between edges 1-2 and 4-5.  In this case, we can replace $c$ and $d$ by classical crossings of opposite sign, so the linking number of any pair of disjoint cycles in $G^*$ is the same as in $G$.  So $C_1$ and $C_2$ must also have odd linking number in $G$, and so $G$ is a virtually linked diagram of $K_6$.  This proves that $K_6$ is in $IVL_2$.  Combined with Figure \ref{F:not3IVL}, we have that $K_6$ is in $IVL_2\backslash IVL_3$.  $\Box$

\begin{prop} \label{P:K331}
$K_{3,3,1}$, $K_{4,4}-edge$ and $G_7$ are in $IVL_2\backslash IVL_3$.
\end{prop}

\noindent{\sc Proof:}  The proofs are similar to Proposition \ref{P:K6}, and the details are left to the reader.  In each case, the proof involves considering a virtual crossing between pairs of non-adjacent edges and showing that, given a diagram with two virtual crossings, either we can find a link with just one virtual crossing, or we can replace the two virtual crossings by classical crossings so as to preserve the sum of the linking numbers in the diagram.  This shows that these graphs are in $IVL_2$, and Figure \ref{F:not3IVL} shows that they are not in $IVL_3$. $\Box$

\begin{figure} [h]
    $$\includegraphics{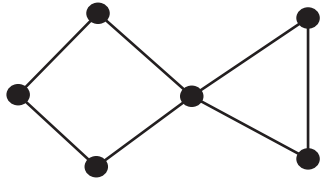}$$
    \caption{The join of a square and a triangle} \label{F:join}
    \end{figure}

We can use the examples of the Petersen graph and $K_6$ to construct graphs which are $n$IVL but not $(n+1)$IVL for all $n$.  Given graphs $G_1$ and $G_2$, the {\it join} of the two graphs, denoted $G_1\wedge G_2$, is the result of choosing a vertex from each graph and identifying them (for example, see Figure~\ref{F:join}).

We will join copies of the Petersen graph $P$ and $K_6$.  First consider $G = \bigwedge_{i=1}^{n}{P}$.  Since $P$ is in $IVL_1\backslash IVL_2$, an unlinked diagram for $G$ will require 2 virtual crossings in each copy of $P$, for a total of $2n$ virtual crossings.  So $G$ is $(2n-1)$IVL but not $2n$IVL.  Now consider $H = \bigwedge_{i=1}^{n-1}{P} \wedge K_6$.  An unlinked diagram for $H$ requires two virtual crossings in each copy of $P$ and {\it three} virtual crossings in the copy of $K_6$, for a total of $2n+1$ virtual crossings.  So $H$ is $2n$IVL but not $(2n+1)$IVL.  Together, these constructions prove the following theorem:

\begin{thm} \label{T:nonempty}
For all $n \geq 1$, the set $IVL_n\backslash IVL_{n+1}$ is not empty.  So the filtration of the sets $IVL_n$ is strictly decreasing and never terminates.
\end{thm}

\section{Virtual Unknotting Number} \label{S:unknotting}

In this section we consider the operation of {\it virtualizing} classical crossings in classical or virtual knots and links, by which we mean replacing the classical crossing by a virtual crossing and leaving the remainder of the diagram unchanged.  Since a purely virtual knot or link is trivial \cite{ka}, virtualizing crossings is an unknotting operation.  This naturally leads us to define the {\it virtual unknotting number}.

\begin{defn}
Given a virtual knot K and a diagram D, define $vu_D(K)$ to be the minimum number of classical crossings in D which need to be virtualized in order to unknot K.  The {\bf virtual unknotting number vu(K)} is the minimum of $\{vu_D(K)\}$, taken over all diagrams D for K.
\end{defn}

Clearly, the virtual unknotting number is no greater than the crossing number.  In fact, it is strictly less than the crossing number, since a virtual knot with exactly one classical crossing is always an unknot.  Like the classical unknotting number, the virtual unknotting number is generally very difficult to compute.  The main result of this section will be to show that any classical knot with non-trivial Jones polynomial has virtual unknotting number at least 2.  This will involve studying how virtualizing a classical crossing changes the Jones polynomial.

\subsection{Virtualizing crossings and the Jones polynomial} \label{SS:jones}

Kauffman \cite{ka} showed how to extend the bracket polynomial, and hence the Jones polynomial, to virtual knots.

\begin{figure} [h]
    $$\includegraphics{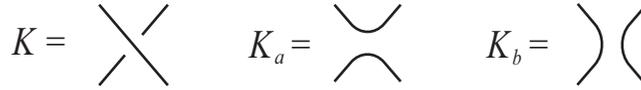}$$
    \caption{Bracket smoothings} \label{F:bracket}
    \end{figure}

\begin{defn}
The {\it bracket polynomial} $\langle K \rangle$ for an unoriented virtual knot K is determined by the following relations:
\begin{itemize}
	\item  $\langle K \rangle = A\langle K_a \rangle + A^{-1}\langle K_b \rangle$
	\item $\langle U_n \rangle = (-A^2-A^{-2})^{n-1}$
\end{itemize}
where $K$, $K_a$ and $K_b$ are as shown in Figure~\ref{F:bracket}, and $U_n$ is the unlink of $n$ components.
\end{defn}

The {\it writhe w(K)} of an oriented virtual knot is just the sum of the signs of the classical crossings.  The Jones polynomial $V_K(t)$ is then defined by:
$$V_K(t) = \left( {(-A^3)^{-w(K)}\langle K \rangle} \right)_{t^{1/2}=A^{-2}}$$
In general, $V_K(t) \in \mathbb{Z}[t^{-1/2},t^{1/2}]$.  However, in classical links with odd numbers of components, including knots, the Jones polynomial has only integer powers of $t$ \cite{li}.

\begin{figure} [h]
    $$\includegraphics{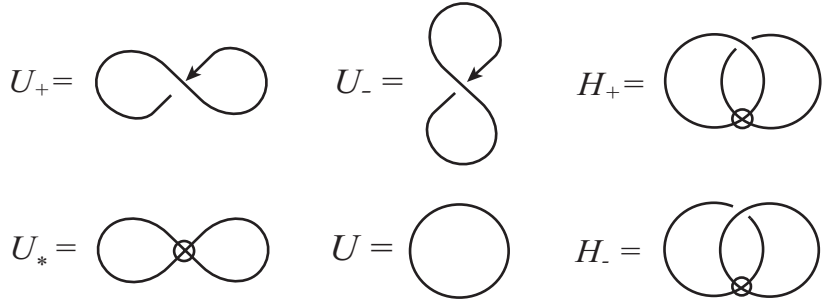}$$
    \caption{Unknots $U_+$, $U_-$, $U_*$ and $U$, and virtual Hopf links $H_+$ and $H_-$} \label{F:unknots}
    \end{figure}

Now we will consider how this polynomial changes when we change the sign of a crossing, or virtualize it.  It will be useful to also consider the Jones polynomial in the form $f_K(A) = (-A^3)^{-w(K)}\langle K \rangle$ (i.e. using the variable $A$ instead of $t$).  We will also want to consider the bracket polynomials of the unknots $U_+$, $U_-$ and $U_*$ and the virtual links $H_+$ and $H_-$ shown in Figure~\ref{F:unknots}.  So $U_+$ is the unknot with a positive twist, $U_-$ is the unknot with a negative twist, and $U_*$ is the unknot with a virtual twist, while $H_+$ and $H_-$ are the two virtual Hopf links.

It is an easy exercise to compute $\langle U_+ \rangle = -A^3$, $\langle U_- \rangle = -A^{-3}$, $\langle U_* \rangle = \langle U \rangle = 1$, and $\langle H_+ \rangle = \langle H_- \rangle = A + A^{-1}$.  Let $K_+$ be a virtual knot with a positive crossing $c$, $K_-$ be the result of changing $c$ to a negative crossing, and $K_*$ be the result of virtualizing $c$.  Now compute $\langle K_+ \rangle$, applying the bracket skein relation to every crossing {\it except} $c$.  Each state in the resulting state sum consists of some number of unlinked circles together with $H_+$, $U_+$ or $U_-$ (the crossing $c$ may have changed its sign, depending on how the arcs were redrawn).  If we collect the terms involving $H_+$, those involving $U_+$ and those involving $U_-$, we see that there will be some polynomials $r(A)$, $s(A)$ and $z(A)$ such that:
$$\langle K_+ \rangle = r(A)\langle H_+ \rangle + s(A)\langle U_+ \rangle + z(A)\langle U_- \rangle = (A + A^{-1})r(A) - A^3s(A) - A^{-3}z(A)$$

Applying the same procedure to $K_-$ and $K_*$, we get the same results, except that in $K_-$ the roles of $U_+$ and $U_-$ are switched, and $H_+$ is replaced by $H_-$, and in $K_*$, $U_+$ and $U_-$ are both replaced by $U_*$, and $H_+$ is replaced by the unlink $U_2$.  So,
$$\langle K_- \rangle = r(A)\langle H_- \rangle + s(A)\langle U_- \rangle + z(A)\langle U_+ \rangle = (A + A^{-1})r(A) - A^{-3}s(A) - A^3z(A)$$
$$\langle K_* \rangle = r(A)\langle U_2 \rangle + s(A)\langle U_* \rangle + z(A)\langle U_* \rangle = (-A^2 - A^{-2})r(A) + s(A) + z(A)$$

Now, let $w$ be the writhe of $K_+$, so the writhes of $K_*$ and $K_-$ are $w-1$ and $w-2$, respectively.  Then:
\begin{align*}
	f_{K_+}(A) &= (-A)^{-3w}((A + A^{-1})r(A) - A^3s(A) - A^{-3}z(A)) \\
			&= -((-A)^{-3w+1} + (-A)^{-3w-1})r(A) + (-A)^{-3w+3}s(A) + (-A)^{-3w-3}z(A) \\
			&= (1+A^{-2}) (-(-A)^{-3w+1}r(A)) + ((-A)^{-3w+3}s(A)) + ((-A)^{-3w-3}z(A))\\
	f_{K_*}(A) &= (-A)^{-3w+3}((-A^2 - A^{-2})r(A) + s(A) + z(A)) \\
			&= -((-A)^{-3w+5} + (-A)^{-3w+1})r(A) + (-A)^{-3w+3}s(A) + (-A)^{-3w+3}z(A) \\
			&= (A^4 + 1)(-(-A)^{-3w+1}r(A)) + ((-A)^{-3w+3}s(A)) + A^6((-A)^{-3w-3}z(A))\\
	f_{K_-}(A) &= (-A)^{-3w+6}((A + A^{-1})r(A) - A^{-3}s(A) - A^3z(A)) \\
			&= -((-A)^{-3w+7} + (-A)^{-3w+5})r(A) + (-A)^{-3w+3}s(A) + (-A)^{-3w+9}z(A) \\
			&= (A^6 + A^4)(-(-A)^{-3w+1}r(A)) + ((-A)^{-3w+3}s(A)) + A^{12}((-A)^{-3w-3}z(A))
\end{align*}

So we have proven the following result (where $t^{1/2}=A^{-2}$):

\begin{thm} \label{T:jones}
Given oriented virtual knots $K_+$, $K_-$ and $K_*$ which differ only in a single crossing (which is positive, negative and virtual in the three knots), there are polynomials $n(t) = -(-t^{-\frac{1}{4}})^{-3w+1}r(t^{-\frac{1}{4}})$, $p(t) = (-t^{-\frac{1}{4}})^{-3w+3}s(t^{-\frac{1}{4}})$ and $q(t) = (-t^{-\frac{1}{4}})^{-3w-3}z(t^{-\frac{1}{4}})$ in $\mathbb{Z}[t^{-1/2},t^{1/2}]$ such that:
$$V_{K_+}(t) = (1+t^{1/2})n(t) + p(t) + q(t)$$
$$V_{K_*}(t) = (t^{-1} + 1)n(t) + p(t) + t^{-3/2}q(t)$$
$$V_{K_-}(t) = (t^{-3/2} + t^{-1})n(t) + p(t) + t^{-3}q(t)$$
\end{thm}

Of particular interest for us is the case when $K_+$ and $K_-$ are {\it classical} knots, so $K_*$ has only one virtual crossing.  In this case, the state sum for $K_+$ does {\it not} involve $H_+$, so $r(A) = 0$, and hence $n(t) = 0$ in Theorem \ref{T:jones}.  So we get the following:

$$V_{K_+}(t) = p(t) + q(t)$$
$$V_{K_*}(t) = p(t) + t^{-3/2}q(t)$$
$$V_{K_-}(t) = p(t) + t^{-3}q(t)$$

This allows us to write down various skein relations for the Jones polynomial for knots with a single virtual crossing.

\begin{cor} \label{C:jones}
Consider the oriented knots (or links) $K_+$, $K_-$, $K_*$ and $K_0$ which differ only at one crossing, which is positive in $K_+$, negative in $K_-$, virtual in $K_*$, and replaced by the oriented smoothing in $K_0$, and with all other crossings classical.  Then we have the following relations among the Jones polynomials of these links:
\begin{itemize}
	\item  $V_{K_+} = (t^{3/2}+1)V_{K_*} - t^{3/2}V_{K_-}$
	\item  $V_{K_-} = (t^{-3/2}+1)V_{K_*} - t^{-3/2}V_{K_+}$
	\item  $V_{K_*} = \frac{V_{K_+} + t^{3/2}V_{K_-}}{t^{3/2}+1}$
	\item  $(t^{-1}+t^{-1/2})V_{K_+} - (t+t^{-1/2})V_{K_*} = (t^{1/2}-t^{-1/2})V_{K_0}$
	\item  $(t^{-1}+t^{1/2})V_{K_*} - (t+t^{1/2})V_{K_-} = (t^{1/2}-t^{-1/2})V_{K_0}$
\end{itemize}
\end{cor}

\noindent{\sc Proof:}  From Theorem \ref{T:jones}, and remembering that $n(t) = 0$, we have that:
$$V_{K_+} - V_{K_-} = (1-t^{-3})q(t) = (1+t^{-3/2})(V_{K_+} - V_{K_*})$$
The first three relationships come from solving this equation for $V_{K_+}$, $V_{K_-}$ and $V_{K_*}$.  The last two relationships come from combining the first three with the well-known skein relation for the Jones polynomial: $t^{-1}V_{K_+} - tV_{K_-} = (t^{1/2}-t^{-1/2})V_{K_0}$.  $\Box$ 

Our next corollary uses the following fact:

\begin{lem} \label{L:integer}
If $K_+$ and $K_-$ are classical knots which differ in a single crossing, and $p(t)$ and $q(t)$ are polynomials such that $V_{K_+}(t) = p(t) + q(t)$ and $V_{K_-}(t) = p(t) + t^{-3}q(t)$ (as in Theorem \ref{T:jones}), then $p(t),\ q(t) \in \mathbb{Z}[t^{-1},t]$ (so all powers of $t$ in $p(t)$ and $q(t)$ are integers).
\end{lem}

\noindent{\sc Proof:}  Since $K_+$ and $K_-$ are classical knots, we know that their Jones polynomials have only integral powers of $t$ \cite{li}.  Say that $p(t) = p_1(t) + t^{1/2}p_2(t)$ and $q(t) = q_1(t) + t^{1/2}q_2(t)$, where $p_i(t),\ q_i(t) \in \mathbb{Z}[t^{-1},t]$.  Then $t^{1/2}p_2(t) + t^{1/2}q_2(t) = 0$ and $t^{1/2}p_2(t) + t^{1/2}t^{-3}q_2(t) = 0$ (so that $V_{K_+}$ and $V_{K_-}$ are left with only integer powers of $t$).  This means that $p_2(t) = -q_2(t) = -t^{-3}q_2(t)$.  But this is only possible if $q_2(t) = 0$, and hence $p_2(t) = 0$.  We conclude that $p(t),\ q(t) \in \mathbb{Z}[t^{-1},t]$.  $\Box$

\begin{cor} \label{C:non-classical}
If $K_+$ and $K_-$ are classical knots which differ in a single crossing, $K_*$ is the result of virtualizing that crossing, and $V_{K_+}(t) \neq V_{K_-}(t)$, then $K_*$ is a non-classical (and hence non-trivial) virtual knot.
\end{cor}

\noindent{\sc Proof:}  From Theorem \ref{T:jones}, setting $n(t) = 0$, we have that $V_{K_+}(t) = p(t) + q(t)$, $V_{K_-}(t) = p(t) + t^{-3}q(t)$ and $V_{K_*}(t) = p(t) + t^{-3/2}q(t)$.  By Lemma \ref{L:integer}, $p(t),\ q(t) \in \mathbb{Z}[t^{-1},t]$.  Moreover, since $V_{K_+}(t) \neq V_{K_-}(t)$, $q(t) \neq 0$.  Therefore, $V_{K_*}(t) = p(t) + t^{-3/2}q(t) \notin \mathbb{Z}[t^{-1},t]$, which implies that $K_*$ is not a classical knot.  $\Box$\\
\\
\noindent{\sc Remark:}  Dye and Kauffman \cite{dk} have used similar techniques to look at the effect on the Jones polynomial of a different way of ``virtualizing" a crossing.

\subsection{Results on the virtual unknotting number} \label{SS:unknotting}

Our main result in this section is the following theorem:

\begin{thm} \label{T:vu1}
If $K$ is a classical knot with $vu(K) = 1$, then $V_K(t) = 1$.
\end{thm}

\noindent{\sc Proof:}  We proceed by contradiction.  Assume $vu(K) = 1$, but $V_K(t) \neq 1$.  Since $K$ is a knot, the powers of $t$ in $V_K(t)$ are all integers.  Since $vu(K) = 1$, we can virtualize some crossing $c$ to get a trivial knot $K_*$, so $V_{K_*}(t) = 1$.  We will consider the case when $c$ is a positive crossing in $K$ (the case when $c$ is negative is similar).

By Theorem \ref{T:jones}, with $n(t) = 0$, $V_K - V_{K_*} = (1-t^{-3/2})q(t) = V_K - 1$, so $t^{-3/2}q(t) = q(t) + 1 - V_K$.  Since $V_K \neq 1$, $q(t) \neq 0$.  By Lemma \ref{L:integer}, $q(t)$ contains only integer powers of $t$.  But then $q(t) + 1 - V_K$ contains only integer powers of $t$, while $t^{-3/2}q(t)$ contains {\it no} integer powers of $t$, so these cannot be equal.  This is the desired contradiction. $\Box$

\begin{conj} \label{C:vu}
If $K$ is a non-trivial classical knot, then $vu(K) > 1$.
\end{conj}

A counterexample to this conjecture would also be an example of a non-trivial knot with trivial Jones polynomial.

\begin{figure} [h]
    $$\includegraphics{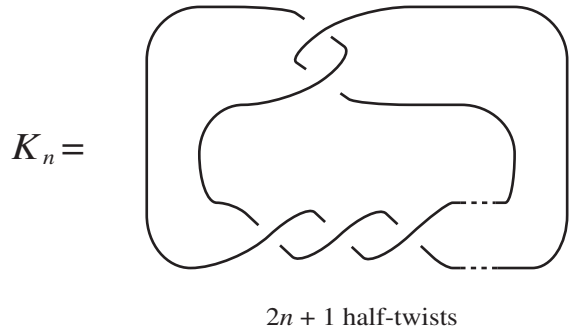}$$
    \caption{The twist knot $K_n$} \label{F:twistknot}
    \end{figure}

As an application of Theorem \ref{T:vu1}, consider the twist knots $K_n$, shown in Figure~\ref{F:twistknot}.  This class of knots includes the trefoil knot as $K_0$.
Using the skein relation $t^{-1}V_{K_+} - tV_{K_-} = (t^{1/2}-t^{-1/2})V_{K_0}$, a straightforward induction shows that $V_{K_n} \neq 1$, so the virtual unknotting number of a twist knot is at least 2.  On the other hand, it is clear that virtualizing the two crossings in the clasp will trivialize $K_n$, so the virtual unknotting number must be exactly 2.

\begin{thm} \label{T:twist}
The virtual unknotting number of the twist knot $K_n$ is 2.
\end{thm}

\section{Intrinsically Virtually Knotted Graphs} \label{S:knotting}

The idea of an {\it intrinsically virtually knotted} graph is completely analogous to the idea of an intrinsically virtually linked graph in Section \ref{S:linking}.

\begin{defn} \label{D:nIVK}
A graph G is {\bf intrinsically virtually knotted of degree n (nIVK)} if every virtual diagram of G with {\bf at most} n virtual crossings contains a virtually knotted cycle.
\end{defn}

As before, we denote the set of $nIVK$ graphs by $IVK_n$.  So $IVK_0$ is the set of (classically) intrinsically knotted graphs, and $IVK_{n+1} \subset IVK_n$, so we have a filtration of sets.  In general, intrinsic knotting is much more difficult to work with than intrinsic linking, and even $IVK_0$ is not yet fully understood.  Our results in this section are not as strong as in Section \ref{S:linking}, and even those we have require the machinery developed in Section \ref{S:unknotting}.  In particular, Theorem \ref{T:vu1} allows us to prove:

\begin{thm} \label{T:ivk1}
If G is an intrinsically knotted graph such that every diagram of G contains a knot with non-trivial Jones polynomial, then G is intrinsically virtually knotted of degree 1.
\end{thm}

\noindent{\sc Proof:}  Consider a diagram $D$ of $G$ with exactly one virtual crossing.  If we replace the virtual crossing with a positive crossing, we get a classical diagram $D_+$ of $G$.  So $D_+$ contains a knotted cycle $C$ with non-trivial Jones polynomial.  By Theorem \ref{T:vu1}, the virtual unknotting number of $C$ is greater than 1.

If the virtual crossing in $D$ is not in the cycle $C$, then $D$ contains $C$ as a knotted cycle.  On the other hand, if the virtual crossing {\it is} in $C$, then the resulting virtually knotted cycle is still non-trivial, since $vu(C) > 1$.  In either case, $D$ contains a virtually knotted cycle, and we conclude that $G$ is intrinsically virtually knotted of degree 1.  $\Box$

It is conjectured that the Jones polynomial distinguishes the unknot, which leads us to conjecture:

\begin{conj} \label{C:ivk0=ivk1}
If a graph G is intrinsically knotted, then G is intrinsically virtually knotted of degree 1.
\end{conj}

While we cannot prove this, we can verify that it holds for all known examples of intrinsically knotted graphs.

\begin{thm} \label{T:K7}
The complete graph on 7 vertices, $K_7$, is intrinsically virtually knotted of degree 1.
\end{thm}

\noindent{\sc Proof:}  Conway and Gordon \cite{cg} showed that every classical diagram of $K_7$ contains a knotted Hamiltonian cycle with Arf invariant 1, so $K_7 \in IVK_0$.  However, there is a well-known connection between the Arf invariant of a knot $K$ and its Jones polynomial $V_K$:  $V_K(\sqrt{-1}) = (-1)^{Arf(K)}$ \cite{li}.  Since ${\rm Arf}(C) = 1$, $V_C$ must be non-trivial.  So, by Theorem \ref{T:ivk1}, $K_7$ is intrinsically virtually knotted of degree 1.  $\Box$

The same argument shows that $K_{3,3,1,1}$ and Foisy's graph $H$ are also intrinsically virtually knotted of degree 1, since the intrinsic knottedness of both graphs was proved by Foisy using the Arf invariant \cite{fo1, fo2}.  The other known intrinsically knotted graphs are derived from $K_7$ and $K_{3,3,1,1}$ using triangle-Y exchanges, in which a triangle (a cycle in the graph of length 3) is replaced by a Y (the three edges of the triangle are removed, and a new vertex is added adjacent to the three vertices of the triangle).  Motwani et. al. \cite{mo} showed that these exchanges preserve intrinsic linking, intrinsic knotting, and many other graph properties.  Their proof easily generalizes to show that if every diagram of a graph $G$ contains a knotted cycle with non-trivial Jones polynomial, and $G'$ is the result of a triangle-Y exchange, then every diagram of $G'$ will also contain a knotted cycle with non-trivial Jones polynomial.  So every graph which can be derived from $K_{3,3,1,1}$ and $K_7$ by triangle-Y exchanges is also intrinsically virtually linked of degree 1.  Since these are all the known minor-minimal intrinsically knotted graphs, we have verified Conjecture \ref{C:ivk0=ivk1} for all known instrinsically knotted graphs.

\begin{figure} [h]
    $$\includegraphics{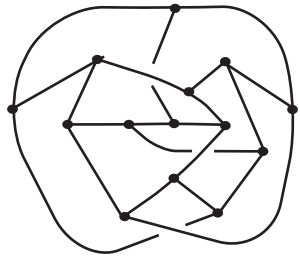}$$
    \caption{The graph $C_{14}$} \label{F:C14}
    \end{figure}

In particular, this means the graph $C_{14}$ shown in Figure~\ref{F:C14} is intrinsically virtually knotted of degree 1, since it is one of the 13 graphs which can be obtained from $K_7$ by triangle-Y exchanges (it is one of the two such graphs which are triangle-free) \cite{ks}.

Since the diagram of $C_{14}$ shown in Figure \ref{F:C14} has only three crossings, virtualizing two of them will mean that the diagram of any cycle contains at most one classical crossing, and is therefore unknotted (the cycle is virtually isotopic to either $U_+$ or $U_-$ of Figure \ref{F:unknots}).  So $C_{14}$ is {\it not} intrinsically virtually knotted of degree 2.  This allows us to prove that the filtration of intrinsically virtually knotted graphs never terminates, similarly to the proof of Theorem \ref{T:nonempty}.

\begin{thm} \label{T:ivknonempty}
For all $n \geq 1$, the set $IVK_{2n-1}\backslash IVK_{2n}$ is not empty.
\end{thm}

\noindent{\sc Proof:}  Let $G = \bigwedge_{i=1}^{n}{C_{14}}$.  Then an unknotted diagram of $G$ requires two virtual crossings in each copy of $C_{14}$.  We conclude that $G \in IVK_{2n-1}\backslash IVK_{2n}$.  $\Box$

\begin{figure} [h]
    $$\includegraphics{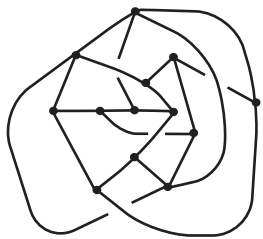}$$
    \caption{The graph $C_{13}$} \label{F:C13}
    \end{figure}

It is still open whether $IVK_{2n}\backslash IVK_{2n+1}$ is likewise non-empty, though it seems likely.  We would need to find a graph which is intrinsically virtually knotted of degree 2, but not of degree 3 (analogous to $K_6$ in the case of intrinsic virtual linking).  One candidate may be the graph $C_{13}$ shown in Figure~\ref{F:C13}.

Since $C_{13}$ is derived from $K_7$ by triangle-Y exchanges, it is $1IVK$; but since the diagram shown in Figure \ref{F:C13} has only four crossings, it is {\it not} $3IVK$.  It is not known, however, whether it is $2IVK$.

\small

\normalsize


\begin{thebibliography}{10}
\bibitem{cg} J. Conway and C. Gordon:  Knots and links in spatial graphs, {\it J. of Graph Theory}, v. 7, no. 4, 1983, pp. 445--453.
\bibitem{dk} H. Dye and L. Kauffman: Minimal surface representations of virtual knots and links, {\it Algebraic and Geometric Topology}, v. 5, 2005, pp. 509-535
\bibitem{ffnp}  E. Flapan, J. Foisy, R. Naimi and J. Pommersheim:  Intrinsically $n$-linked graphs, {\it J. Knot Theory and its Ramifications}, v. 10, no. 8, 2001, pp. 1143-1154
\bibitem{fm}  T. Fleming and B. Mellor:  Virtual Spatial Graphs, preprint, 2005, {\it arXiv:math.GT/0510158}
\bibitem{fo1} J. Foisy: Intrinsically Knotted Graphs, {\it J. Graph Theory}, v. 39, no. 3, 2002, pp. 178-187
\bibitem{fo2} J. Foisy:  A Newly Recognized Intrinsically Knotted Graph, {\it J. Graph Theory}, v. 43, no. 3, 2003, pp. 199-209
\bibitem{gu} R. Guy: The decline and fall of Zarankiewicz's theorem, in {\it Proof Techniques in Graph Theory} (F. Harary Ed.), Academic Press, New York, 1969, pp. 63-69.
\bibitem{ka}  L. Kauffman:  Virtual Knot Theory, {\it Europ. J. Combinatorics}, v. 20, 1999, pp. 663-691
\bibitem{ka2}  L. Kauffman:  Invariants of Graphs in Three-Space, {\it Trans. Amer. Math. Soc.}, v. 311, no. 2, 1989, pp. 697-710
\bibitem{ks} T. Kohara and S. Suzuki:  Some remarks on knots and links in spatial graphs, {\it Knots 90}, de Gruyter \& Co., Berlin, 1992, pp. 435-445
\bibitem{li}  W.B.R. Lickorish, {\it An Introduction to Knot Theory}, Springer, New York, 1997
\bibitem{mo}  R. Motwani, A. Raghunathan, and H. Saran: Constructive Results from Graph Minors: Linkless Embeddings, {\it 29th Annual Symposium on Foundations of Computer Science}, IEEE, 1988, pp. 398-409
\bibitem{rst} ÊN. Robertson, P. Seymour, and R.Thomas: Sachs' linkless embedding conjecture, {\it  J. Combin. Theory Ser. B} {\bf 64} (1995), no. 2, 185--227. 
\bibitem{ya}  S. Yamada:  An Invariant of Spatial Graphs, {\it J. Graph Theory}, v. 13, no. 5, 1989, pp. 537-551
\end{thebibliography}
\end{document}